\documentstyle[12pt, amssymb, amsmath]{article}
\begin{document}
\newcommand{\bref}[1]{$\mbox{(\ref{#1})}$}
\newcommand{\s}{\sqrt {n + b}}
\newcommand{\w}{\sqrt {t + b}}
\newcommand{\p}{{\Bbb {R^+}} \rightarrow {\Bbb {C}}}
\newcommand{\n}{\sqrt {2n + 1}}
\newcommand{\r}{{\Bbb {R}} ^d}
\newcommand{\q}{{\mathcal Q}}
\newcommand{\z}{{\Bbb {Z}} ^d}

\begin{center}
\bf On the $L^2$-norm of periodizations of functions

~

\rm Oleg Kovrijkine\\
School of Mathematics\\
Institute for Advanced Study\\
Princeton, NJ 08540, USA\\
{\sl E-mail address}: olegk@ias.edu\\

~

To the memory of Tom Wolff\\
\end{center}

\begin{abstract} We prove that the $L^2([0,1]^d \times SO(d))$-norm of
periodizations
 of a function from $L^1({\Bbb {R}} ^d)$ is equivalent to the $L^2({\Bbb {R}}
^d)$-norm of the function itself in higher dimensions. We generalize the statement
 for functions from $L^p({\Bbb {R}} ^d)$ where $1 \le p < \frac {2d}{d + 2}$
 in the spirit of the Stein-Tomas theorem.
\end{abstract}

{\bf 0. Introduction.}\\

Let $f$ be a function from $L^1({\Bbb {R}} ^d)$. Define a family of its
 periodizations with respect to a rotated integer lattice:
\begin{eqnarray} g_{\rho}(x) = \sum \limits _{\nu \in {\Bbb {Z}} ^d} f(\rho (x -
\nu))\label{g}\end{eqnarray}
for all rotations $\rho \in SO(d)$. The main object of our study is $G$,
 the $L^2([0,1]^d \times SO(d))$-norm of the family of periodizations,
\begin{eqnarray} G^2 &=& \int\limits_{\rho \in
SO(d)}\int\limits_{[0,1]^d}|g_{\rho}(x)|^2 dx d\rho\nonumber \\
&=& \int\limits_{\rho \in SO(d)}\|g_{\rho}\|_2^2d\rho \label{g10}.\end{eqnarray}

    The purpose of this work is to show how $G$ can give an estimate of
 the $L^2({\Bbb {R}} ^d)$-norm of a function from $L^1({\Bbb {R}} ^d)$
 in higher dimensions. Some results on the Steinhaus tiling problem
 are related to {\bf Theorem 1} since periodizations naturally appear
 in the problem of Steinhaus. M. Kolountzakis (\cite{K}) proves that if
 a function $f \in L^1({\Bbb {R}} ^2)$ and $|x|^{\alpha}f \in L^1({\Bbb {R}}
^2)$,
 where $\alpha > \frac {10} 3$ and its periodizations are constants, then the
function is continuous. Another result is obtained by M. Kolountzakis and T.
Wolff (\cite{KW}, Theorem 1).
 It  says that if periodizations of a function from $L^1({\Bbb {R}} ^d)$ are
constants then the function is continuous provided that the dimension d is at
least three.  \\

The main theorems are the following.

~

{\bf Theorem 1}:
{\it let $d \ge 4$ and let $f \in L^1({\Bbb {R}} ^d)$. If periodizations of f
\begin{eqnarray*} g_{\rho}(x) = \sum \limits _{\nu \in {\Bbb {Z}} ^d} f(\rho (x
- \nu))\end{eqnarray*}
are in $L^2([0,1]^d)$ for almost all rotations $\rho \in SO(d)$ and
\begin{eqnarray*} G^2 = \int\limits_{\rho \in SO(d)}\|g_{\rho}\|_2^2d\rho <
 \infty \end{eqnarray*}
then $f \in L^2({\Bbb {R}} ^d)$:
$$\|f\|_2 \le C(G + \|f\|_1),$$
where $C$ depends only on $d$.}

~

We also obtain the following inverse theorem.

~

{\bf Theorem 1}$'$: {\it let $d \ge 5$, let $f \in L^1({\Bbb {R}} ^d) \cap L^2({\Bbb {R}}
^d)$, and
 let $g_{\rho}$ be periodizations of f
\begin{eqnarray*} g_{\rho}(x) = \sum \limits _{\nu \in {\Bbb {Z}} ^d} f(\rho (x
- \nu))\end{eqnarray*}
then $g_{\rho} \in L^2([0,1]^d)$ for almost all rotations $\rho \in SO(d)$ and
\begin{eqnarray}\int\limits_{\rho \in SO(d)}\|g_{\rho}\|_2^2 d\rho  \le
C(\|f\|_2 + \|f\|_1)^2, \label{G88}\end{eqnarray} where $C$ depends only on $d$.}

~

We will generalize {\bf Theorems 1} and {\bf 1}$'$ in the spirit of the
Stein-Tomas Theorem (\cite{DC}, Chapter 6.5).

~

{\bf Theorem 2}: {\it let $d \ge 4$, let $1 \le p < \frac {2d} {d + 2}$, and let $f \in
L^1({\Bbb {R}} ^d) \cap L^p({\Bbb {R}} ^d)$. If periodizations of f
\begin{eqnarray*} g_{\rho}(x) = \sum \limits _{\nu \in {\Bbb {Z}} ^d} f(\rho (x
- \nu))\end{eqnarray*}
are in $L^2([0,1]^d)$ for almost all rotations $\rho \in SO(d)$ and
\begin{eqnarray*} G^2 = \int\limits_{\rho \in SO(d)}\|g_{\rho}\|_2^2d\rho <
 \infty \end{eqnarray*}
then $f \in L^2({\Bbb {R}} ^d)$:
\begin{eqnarray}\|f\|_2 \le C(G + \|f\|_p), \label{G89}\end{eqnarray}
where $C$ depends only on $d$ and $p$.}

~

We also obtain the following inverse theorem.

~

{\bf Theorem 2}$'$: {\it let $d \ge 5$, let $1 \le p < \frac {2d} {d + 2}$, and let $f
\in L^1({\Bbb {R}} ^d) \cap L^2({\Bbb {R}} ^d)$,
  and let $g_{\rho}$ be periodizations of f
\begin{eqnarray*} g_{\rho}(x) = \sum \limits _{\nu \in {\Bbb {Z}} ^d} f(\rho (x
- \nu))\end{eqnarray*}
then $g_{\rho} \in L^2([0,1]^d)$ for almost all rotations $\rho \in SO(d)$ and
\begin{eqnarray}\int\limits_{\rho \in SO(d)}\|g_{\rho} - \hat g_{\rho} (0)\|_2^2
d\rho  \le C(\|f\|_2 + \|f\|_p)^2, \label{G90}\end{eqnarray} where $C$ depends only on
$d$ and $p$.}

~
The rest of the paper is concerned with the proofs of the theorems stated
above.\\

{\bf 1. Case $p = 1$.}\\

    Note that the constant $C$ below is not fixed and varies appropriately
 from one equality or inequality to another although such variations are not noted.

~

Proof of {\bf Theorem 1}:\\

We will denote $\tilde f(x) = \bar f (-x)$ and $F(x) = f * \tilde f (x)$. Then
$F \in L^1(\r)$ and
\begin{eqnarray}\|F\|_1 \le \|f\|^2_1.\label{F}\end{eqnarray}
We will define the following functions $h, h_1, h_2: {\Bbb {R^+}} \rightarrow
{\Bbb {C}}$
\begin{eqnarray}h(t) &=& \int |\hat f(\xi)|^2 d\sigma_t(\xi)\label{h}\\
&=& \int_{\r} f*\tilde f (x) \widehat {d\sigma_t}(x)dx \nonumber \\
&=& \int_{\r} F (x) \widehat {d\sigma_t}(x)dx, \label{h3}\end{eqnarray}
\begin{eqnarray}h_1(t)
= \int_{|x| \le 1} F (x) \widehat {d\sigma_t}(x)dx, \label{h1}\end{eqnarray}
\begin{eqnarray}h_2(t)
= \int_{|x| > 1} F (x) \widehat {d\sigma_t}(x)dx. \label{h2}\end{eqnarray}
Clearly $h = h_1 + h_2$.

~

{\bf Lemma 1} {\it
Let $q: {\Bbb {R}} \rightarrow {\Bbb {R}}$ be a Schwartz function supported in
 $[\frac 12, 2]$, and let $b \in [0,1)$. Define $H_1:{\Bbb {R}} \rightarrow
{\Bbb {C}}$
\begin{eqnarray*} H_1 (t) = \frac 1 {\w} h_1(\w) q(\frac {\w} N).\end{eqnarray*}
Then for large enough $N$ we have

\begin{eqnarray} \sum \limits _{\nu \neq 0} |\hat H_1 (\nu)| \le \frac
{C\|F\|_1} {N} \label{H1}\end{eqnarray}
where $C$ depends only on $q$ and $d$.}

~

Proof of {\bf Lemma 1}:\\

First we will estimate derivatives of $h_1 (t)$
\begin{eqnarray}|h_1^{(k)} (t)| \le Ct^{d-1}\|F\|_1\label{dif}\end{eqnarray}
where $t \ge 1$ and $C$ depends only on $k$ and $d$. This follows from \bref{h1} by
differentiating the last equality $k$ times:
\begin{eqnarray*}h_1 (t) &=& \int_{|x| \le 1} F (x) \widehat {d\sigma_t}(x)dx\\
&=&t^{d-1}\int_{|x| \le 1} F (x)\int _{|\xi| = 1} e^{-i2\pi t x\cdot \xi}d\sigma(\xi)dx.
\end{eqnarray*}

We can easily prove by induction that
\begin{eqnarray}\frac {d^k} {dt^k} \left (\frac {h_1 (\w)} {\w} \right) =
 \sum \limits _{i = 0} ^ {k} C_{i,k}\frac {h_1^{(i)}(\w)}{(\w)^{2k + 1
-i}}.\label {dif1}\end{eqnarray}
It follows from \bref{dif1} and \bref{dif} that when $t \thicksim N^2$ we have
\begin{eqnarray}\left |\frac {d^k} {dt^k} \left (\frac {h_1 (\w)} {\w} \right)
\right|
 \le CN^{d-k-2}\|F\|_1\label{dif3}\end{eqnarray}
with $C$ depending only on $k$ and $d$.\\

Since $q(\frac {(\w)} {N}) = q(\sqrt {t' + b'}) = \tilde q(t')$ with $t' = \frac
t {N^2}$
 and $b' = \frac b {N^2}$ and $\tilde q(t')$ is a Schwartz function supported in
 $t' \thicksim 1$, we have
\begin{eqnarray}\left |\frac {d^k} {dt^k} q(\frac {(\w)} {N}) \right| &=&
N^{-2k}|\frac {d^k} {d{t'}^k} \tilde q (t')|\nonumber \\
&\le& CN^{-2k}\label {dif2}\end{eqnarray}
with $C$ depending only on $k$ and $q$.\\

Since $q(\frac {(\w)} {N})$ is supported in $t \thicksim N^2$ it follows from
\bref{dif3} and \bref{dif2} that
\begin{eqnarray}\left |\frac {d^k} {dt^k} H_1(t)\right| &=&
\left |\frac {d^k} {dt^k} \left (\frac {h_1 (\w)} {\w} q(\frac {\w} {N})\right)
\right|
\nonumber \\
 &\le& CN^{d-2-k}\|F\|_1\label{dif4}\end{eqnarray}
with $C$ depending only on $k$, $d$ and $q$. Since $H_1(t)$ is also supported in
$t \thicksim N^2$ we have
$$\|H_1^{(k)}\|_1 \le CN^{d-k}\|F\|_1.$$
Therefore
\begin{eqnarray} |\hat H_1(\nu)| &\le& \frac {C} {|\nu|
^k}\|H_1^{(k)}\|_1\nonumber \\
&\le& \frac {CN^{d-k}\|F\|_1} {|\nu| ^k}\label{nu}\end{eqnarray}
for every $\nu \neq 0$.\\

Summing \bref{nu} over all $\nu \neq 0$ and putting $k = d + 1$ we get our
desired result
\begin{eqnarray} \sum \limits _{\nu \neq 0} |\hat H_1 (\nu)| \le \frac
{C\|F\|_1} {N}
\end{eqnarray}
where $C$ depends only on $q$ and $d$.\hfill$\square$

~

In the next lemma we will use an approach related to (\cite{KW}, Lemma 1.1).\\
{\bf Lemma 2} {\it
Let $q: {\Bbb {R}} \rightarrow {\Bbb {R}}$ be a Schwartz function supported in
 $[\frac 12, 2]$, and let $b \in [0,1)$. Define $H_2:{\Bbb {R}} \rightarrow
{\Bbb {C}}$
\begin{eqnarray*} H_2 (t) = \frac 1 {\w} h_2(\w) q(\frac {\w} N).\end{eqnarray*}
Then for large enough $N$ we have

\begin{eqnarray} \sum \limits _{\nu \neq 0} |\hat H_2 (\nu)| \le \int
\limits_{|x| \ge 1}|F(x)|\cdot |D_N(x)|\label{H2}\end{eqnarray}
where $D_N: \r \rightarrow \Bbb {C}$
\begin{eqnarray}|D_N(x)| \le C \begin{cases} (\frac N {|x|})^{\frac {d-2}
2}&\text{if $|x| \ge \frac N 2$}\\
\frac 1 {N}&\text{if $1 \le |x| \le \frac N 2$,}
\end{cases} \label{D}\end{eqnarray} with $C$ depending only on $q$ and $d$.}

~

Proof of {\bf Lemma 2}:\\

We have
\begin{eqnarray}&&\hat H_2 (\nu)\nonumber \\
 &=& \int H_2(t)e^{-i2\pi \nu t}dt\nonumber \\
 &=& 2e^{i2\pi \nu b}\int Nq(t) h_2(tN)e^{-i2\pi \nu (Nt)^2}dt\nonumber \\
&=& 2e^{i2\pi \nu b}\int Nq(t)e^{-i2\pi \nu (Nt)^2}\int_{|x| > 1} F (x) \widehat
{d\sigma_{Nt}}(x)dxdt\nonumber \\
&=& 2e^{i2\pi \nu b}\int_{|x| > 1} F (x)\int Nq(t)e^{-i2\pi \nu (Nt)^2}
(Nt)^{d-1}\widehat {d\sigma}(Ntx)dtdx.
\label{H3} \end{eqnarray}

We will use a well-known fact that $\widehat {d\sigma} (x) = Re(B(|x|))$ with
$B(r)
 = a(r)e^{i2\pi r}$ and $a(r)$ satisfying estimates
\begin{eqnarray}|a^{k}(r)| \le \frac {C}{r^{\frac{d-1}2 +
k}}\label{a}\end{eqnarray}
with $C$ depending only on $k$ and $d$.
Now we will need to estimate the inner integral in \bref{H3} with $B(|x|)$
instead of
 $\widehat {d\sigma}(x)$
\begin{eqnarray}&&\int Nq(t)e^{-i2\pi \nu (Nt)^2}(Nt)^{d-1} a(N|x|t)e^{i2\pi
N|x|t }dt
\nonumber \\ &=&\frac {N^{\frac {d+1} 2}} {|x|^{\frac {d-1} 2}}\int q(t)
e^{-i2\pi \nu (Nt)^2}t^{d-1} a(N|x|t)(N|x|)^{\frac {d-1} 2}e^{i2\pi N|x|t
}dt\nonumber \\
&=&\frac {N^{\frac {d+1} 2}} {|x|^{\frac {d-1} 2}}e^{i2\pi \frac
{|x|^2}{4\nu}}\int q(t) a(N|x|t)(N|x|)^{\frac {d-1} 2}t^{d-1}e^{-i2\pi \nu N^2
(t-\frac {|x|}{2\nu N})^2}dt\nonumber \\
&=&\frac {N^{\frac {d+1} 2}} {|x|^{\frac {d-1} 2}}e^{i2\pi \frac
{|x|^2}{4\nu}}\int \phi(t)e^{-i2\pi \nu N^2 (t-\frac {|x|}{2\nu N})^2}dt
\label{a1}\end{eqnarray}
where $\phi(t) = q(t) a(N|x|t)(N|x|)^{\frac {d-1} 2}t^{d-1}$ is a Schwartz
function
 supported in $[\frac 12, 2]$ whose derivatives and the function itself are
bounded
 uniformly in $t$, $x$ and $N$ because of \bref{a}. Note that we used here the
fact
 that $N|x| \ge 1$. We can say even more. Note that in fact $\phi(t) = \phi(t,
|x|)$.
 Let $|x| = c\cdot r$ where $c \ge 2$ and $r \ge \frac 1 2$. Then all partial
derivatives
 of $\phi(t, c\cdot r)$ with respect to $t$ and $r$ are also bounded uniformly
in $t$,
 $r$, $c$ and $N$. The only place, where we will use that $\phi(t)$ also depends on $x$,
 is formula \bref{psi10} from the proof of {\bf Lemma 4}. Therefore, we will keep writing
just
 $\phi(t)$ until formula \bref{psi10}.\\
 
From the method of stationary phase (\cite{H}, Theorem 7.7.3) it follows that if
$k \ge 1$ then
\begin{eqnarray}|\int \phi(t)e^{-i2\pi \nu N^2 (t-\frac {|x|}{2\nu N})^2}dt -
\sum \limits _{j=0}^{k-1}c_j(\nu N^2)^{-j-\frac 12}\phi^{(2j)}(\frac {|x|}{2\nu
N})|
 \le c_k(\nu N^2)^{-k-\frac 12}\label{phi}\end{eqnarray}
where $c_j$ are some constants.

Since $\phi$ is supported in $[\frac 12, 2]$ we conclude from \bref{phi} that
\begin{eqnarray}|\int \phi(t)e^{-i2\pi \nu N^2 (t-\frac {|x|}{2\nu N})^2}dt| \le
\begin{cases} C(\nu N^2)^{-\frac 12}&\text{if $\nu \in [\frac {|x|} {4N},\frac
{|x|} {N}]$}\\
C_k(\nu N^2)^{-k-\frac 12}&\text{if $\nu \notin [\frac {|x|} {4N},\frac {|x|}
{N}]$}\end{cases}
.\label{phi1}\end{eqnarray}

If $\frac {|x|} {N} \le \frac 12$, then there are no $\nu$ in
$[\frac {|x|} {4N},\frac {|x|} {N}]$ and therefore if we sum \bref{phi1} over
 all $\nu \neq 0$ we will get
\begin{eqnarray}|\int \phi(t)e^{-i2\pi \nu N^2 (t-\frac {|x|}{2\nu N})^2}dt| \le
C_kN^{-2k-1}.\label{phi2}\end{eqnarray}

If $\frac {|x|} {N} \ge \frac 12$ then the number of $\nu$ in $[\frac {|x|}
{4N},
\frac {|x|} {N}]$ is bounded by $\frac {|x|} {N}$ and therefore if we sum
\bref{phi1}
 over all $\nu \neq 0$ we will get
\begin{eqnarray}\sum \limits _{\nu \neq 0}|\int \phi(t)e^{-i2\pi \nu N^2
(t-\frac {|x|}{2\nu N})^2}dt| &\le& C\frac {|x|} {N} (|x|N)^{-\frac 12} +
C_kN^{-2k-1}\nonumber \\
&\le& C_k\frac {|x|^{\frac 12}} {N^{\frac 32}}.\label{phi3}\end{eqnarray}

Summing \bref{a1} over all $\nu \neq 0$ and applying \bref{phi2} or \bref{phi3}
we conclude
\begin{eqnarray}\sum \limits _{\nu \neq 0} |\int Nq(t)e^{-i2\pi \nu
(Nt)^2}(Nt)^{d-1}
 B(N|x|t)dt| \le \begin{cases} C_k(\frac N {|x|})^{\frac {d-2} 2}&\text{if
$\frac {|x|} N
 \ge \frac 12$}\\
C_k\frac {N^{\frac {d+1} 2 -2k -1}}{|x|^{\frac {d-1} 2}}&\text{if $\frac {|x|} N
\le \frac 12$}\end{cases}
.\label{a2}\end{eqnarray}
Replacing in \bref{H3} $\widehat {d\sigma}(x)$ with $\frac {B(|x|) + \bar
B(|x|)} 2$,
 summing over all $\nu \neq 0$ and applying \bref{a2} with $k \ge \frac {d+1} 4$
we
 get the desired result
\begin{eqnarray} \sum \limits _{\nu \neq 0} |\hat H_2 (\nu)| \le \int
\limits_{|x| \ge 1} |F(x)|\cdot |D_N(x)|\end{eqnarray}
where $D_N: \r \rightarrow \Bbb {C}$
\begin{eqnarray}|D_N(x)| \le C\begin{cases} (\frac N {|x|})^{\frac {d-2}
2}&\text{if $|x| \ge \frac N 2$}\\
\frac 1 {N}&\text{if $1 \le |x| \le \frac N 2$}
\end{cases}\end{eqnarray} with $C$ depending only on $q$ and $d$.\hfill$\square$

~

Now we are in a position to proceed with the proof of {\bf Theorem 1}.
From \bref {g} it follows that
\begin{eqnarray} \hat g_{\rho}(m) = \hat f (\rho m)\label{g1}\end{eqnarray}
for every $m \in \z$. By scaling we can assume that
\begin{eqnarray} \hat g_{\rho}(m) = \hat f (\frac {\rho m} {\sqrt
2})\label{g2}\end{eqnarray}
for every $m \in \z$. It follows that
\begin{eqnarray} \|g_{\rho}\|_2^2 = \sum \limits _{m \in \z} |\hat
g_{\rho}(m)|^2 =
 \sum \limits _{m \in \z}|\hat f (\frac {\rho m} {\sqrt 2})|^2.
\label{g111}\end{eqnarray} Let $r_d (n)$ denote the number of representations of an
integer $n$ as sums of $d$ squares.
 It is a well-known fact from Number Theory that if $d \ge 5$ then
\begin{eqnarray} r_d (n) \ge C n^{\frac {d-2}{2}}\label{r}\end{eqnarray}
and if $d = 4$ and $n$ is odd then
\begin{eqnarray} r_4 (n) \ge C n\label{r1}\end{eqnarray}
where $C > 0$ depends only on $d$.
See for example (\cite{G}, p.30, p.155, p.160).\\

Integrating \bref{g111} with respect to the Haar measure $d\rho$ and applying \bref{g10}
we have
\begin{eqnarray} G^2 &=& \int \limits _{\rho \in SO (d)} \sum \limits _{m \in
\z}|\hat f (\frac {\rho m} {\sqrt 2})|^2 d\rho\nonumber \\
&=&\int \limits _{|\xi| = 1} \sum \limits _{m \in \z}|\hat f (\frac {|m|} {\sqrt
2}\xi)|^2 d\sigma(\xi)\nonumber \\
&=&\sum \limits _{n \ge 0}\sum \limits _{|m|^2 = n}\int \limits _{|\xi| = 1}
|\hat f
 (\frac {|m|} {\sqrt 2}\xi)|^2 d\sigma(\xi)\nonumber \\
&\ge&\sum \limits _{n \ge 0}\sum \limits _{|m|^2 = 2n + 1}\int \limits _{|\xi| =
1}
|\hat f (\frac {|m|} {\sqrt 2}\xi)|^2 d\sigma(\xi)\nonumber \\
&=&\sum \limits _{n \ge 0}r_d({2n+1})\int \limits _{|\xi| = 1} |\hat f
(\sqrt{n + \frac 12}\xi)|^2 d\sigma(\xi)\nonumber \\
&=&\sum \limits _{n \ge 0}\frac {r_d({2n+1})}{(n + \frac 12)^{\frac {d-1}2}}\int
|\hat f (\xi)|^2 d\sigma_{\sqrt{n + \frac 12}}(\xi)
. \label{g21}\end{eqnarray}
Using \bref{h} and \bref{r} or \bref{r1} we conclude from \bref{g21} that
\begin{eqnarray}\sum \limits_{n \ge 0} \frac 1 {\sqrt{n + \frac 12}} h(\sqrt{n +
\frac 12})
\le CG^2.\label{h10}\end{eqnarray}

Let $q: {\Bbb {R}} \rightarrow {\Bbb {R}}$ be a fixed non-negative Schwartz
function
supported in $[\frac 12, 2]$ such that
$$q(x) + q(x/2) = 1$$
when $x \in [1,2]$. It follows that
\begin{eqnarray}\sum \limits_{j \ge 0} q(\frac x {2^j}) = 1
\label{q10}\end{eqnarray}
when $x \ge 1$.\\

Applying the Poisson summation formula to
\begin{eqnarray*}H(t) &=& \frac 1 {\sqrt {t + \frac 12}}h(\sqrt {t + \frac
12})q(\frac
{\sqrt {t + \frac 12}} N)\\
&=& \frac 1 {\sqrt {t + \frac 12}}h_1(\sqrt {t + \frac 12})q(\frac {\sqrt {t +
\frac 12}} N)
+ \frac 1 {\sqrt {t + \frac 12}}h_2(\sqrt {t + \frac 12})q(\frac {\sqrt {t +
\frac 12}} N)\\
&=& H_1(t) + H_2(t)
\end{eqnarray*}
we have
\begin{eqnarray}\sum \limits_{n} H(n) &=&
\sum \limits_{\nu} \hat H(\nu)\nonumber \\
&=& \hat H(0) +
\sum \limits_{\nu \neq 0} \hat H_1(\nu) + \sum \limits_{\nu \neq 0} \hat
H_2(\nu).
 \label{h11}\end{eqnarray}
Note that
\begin{eqnarray}\hat H(0) &=& \int \frac 1 {\sqrt {t + \frac 12}}h(\sqrt {t +
\frac 12})q(\frac {\sqrt {t + \frac 12}} N) dt\nonumber \\
 &=& 2\int h(t)q(\frac {t} N) dt.
 \label{h12}\end{eqnarray}
Substituting \bref{h12} into \bref{h11} we get that
\begin{eqnarray}2\int h(t)q(\frac {t} N) dt &\le& \nonumber \\ \sum \limits_{n}
H(n) + \sum \limits_{\nu \neq 0} |\hat H_1(\nu)| + \sum \limits_{\nu \neq 0}
|\hat H_2(\nu)|
&\le& \nonumber \\
 \sum \limits_{n \ge 0}\frac 1 {\sqrt {n + \frac 12}}h(\sqrt {n + \frac
12})q(\frac
{\sqrt {n + \frac 12}} N) + \frac {C\|F\|_1} {N} + \int\limits_{|x| \ge 1}
|F(x)D_N(x)|dx\label{h35} \end{eqnarray}
where the last inequality follows from {\bf Lemma 1} and {\bf Lemma 2}.\\

From the definition of $D_N(x)$ in \bref{D}, it follows that
\begin{eqnarray}\sum\limits_{j \ge 0}|D_{2^j}(x)| &=& \sum\limits_{2^j \le
2|x|}|D_{2^j}(x)| + \sum\limits_{2^j > 2|x|}|D_{2^j}(x)|\nonumber\\
&\le& \sum\limits_{2^j \le 2|x|}C(\frac {2^j}{|x|})^{\frac {d-2}2} +
\sum\limits_{2^j >
 2|x|}\frac C {2^j} \le C
 \label{D2}\end{eqnarray}
for every $|x| \ge 1$.\\

Putting $N=2^j$ in \bref{h35}, summing over all $j \ge 0$ and applying
\bref{q10} we get by Lebesgue Monotone Convergence Theorem
\begin{eqnarray}2\int\limits_1^{\infty} h(t) dt &\le& \sum \limits_{n \ge
0}\frac 1 {\sqrt {n + \frac 12}}h(\sqrt {n + \frac 12}) + C\|F\|_1 +
C\int\limits_{|x| \ge 1} |F(x)|dx \nonumber \\
&\le& C(G^2 + \|F\|_1)
\label{h15} \end{eqnarray}
where the last inequality follows from \bref{h10}. From the definition of $h(t)$
\bref{h} it follows that
\begin{eqnarray}h(t) \le C\|f\|_1^2
\label{h16} \end{eqnarray}
for $t \le 1$. Therefore we have
\begin{eqnarray}\int |\hat f(x)|^2 dx &=& \int\limits_0^{\infty} |\hat f(\xi)|^2
d\sigma_{t}(\xi)dt\nonumber \\
&=& \int\limits_0^{\infty} h(t) dt \nonumber \\
&\le& C(G^2 + \|f\|_1^2)
\label{f10} \end{eqnarray}
where the last inequality is obtained from \bref{h15}, \bref{h16} and \bref{F}.
From \bref{f10} it follows that $f \in L^2$ and
$$\|f\|_2 \le C(G + \|f\|_1)$$
with $C$ depending only on $d$.\\
\hfill$\square$\\

~

If $d \ge 5$ then
\begin{eqnarray} r_d (n) \le C n^{\frac {d-2}{2}}\label{r2}\end{eqnarray}
where $C > 0$ depends only on $d$. See for example (\cite{G}, p.155, p.160). An argument
similar to one used to get \bref{g21},
 but without scaling, shows that
\begin{eqnarray} G^2 &=& \int \limits _{\rho \in SO (d)} \sum \limits _{m \in
\z}|\hat f
 ({\rho m})|^2 d\rho\nonumber \\
&=&|\hat f (0)|^2 + \sum \limits _{n \ge 1}r_d({n})\int \limits _{|\xi| = 1}
|\hat f (\sqrt{n}\xi)|^2 d\sigma(\xi)\nonumber \\
&=&|\hat f (0)|^2 + \sum \limits _{n \ge 1}\frac {r_d({n})}{n^{\frac
{d-1}2}}\int
|\hat f (\xi)|^2 d\sigma_{\sqrt{n}}(\xi)
. \label{g20}\end{eqnarray}
Using \bref{h} and \bref{r2} we conclude from \bref{g20} that
\begin{eqnarray}G^2 \le \|f\|_1^2 + C\sum \limits_{n \ge 1} \frac 1 {\sqrt{n}}
h(\sqrt{n}).\label{h20}\end{eqnarray}
Repeating arguments which we used to obtain \bref{h15} we get
\begin{eqnarray}\sum \limits_{n \ge 1} \frac 1 {\sqrt{n}} h(\sqrt{n}) &\le&
2\int\limits_{0}^{\infty}h(t)dt + C\|F\|_1\nonumber\\
&\le& C(\|\hat f\|_2^2 + \|f\|_1^2).\label{h22}\end{eqnarray}
Hence we can formulate an inverse theorem to {\bf Theorem 1}:

~

{\bf Theorem 1}$'$:
{\it let $d \ge 5$ and let $f \in L^1({\Bbb {R}} ^d) \cap L^2({\Bbb {R}} ^d)$
and
let $g_{\rho}$ be periodizations of f
\begin{eqnarray} g_{\rho}(x) = \sum \limits _{\nu \in {\Bbb {Z}} ^d} f(\rho (x -
\nu))\label{g30}\end{eqnarray}
then $g_{\rho} \in L^2([0,1]^d)$ for almost all rotations $\rho \in SO(d)$ and
$$\int\limits_{\rho \in SO(d)}\|g_{\rho}\|_2^2 d\rho  \le C(\|f\|_2 +
\|f\|_1)^2$$
where $C$ depends only on $d$.}

~

{\bf Corollary}: complex interpolation between the trivial $p=1$ and $p=2$ gives
us
the following result for $1 < p < 2$: let
$d \ge 5$ and let
$f \in L^{\frac p 2}({\Bbb {R}} ^d) \cap L^p({\Bbb {R}} ^d)$ and let $g_{\rho}$
be
periodizations of f
\begin{eqnarray*} g_{\rho}(x) = \sum \limits _{\nu \in {\Bbb {Z}} ^d} f(\rho (x
- \nu))\end{eqnarray*}
then $g_{\rho} \in L^p([0,1]^d)$ for almost all rotations $\rho \in SO(d)$ and
$$\left ( \int\limits_{\rho \in SO(d)}\|g_{\rho}\|_p^{p'} d\rho \right )^{1/p'}
\le C\|f\|_p^{2-p}(\|f\|_p +
\|f\|_{\frac p 2})^{p-1} = C \|f\|_p(1 +
\frac {\|f\|_{\frac p 2}} {\|f\|_p})^{p-1}$$
where $C$ depends only on $d$.
\\

If $p'$ is an even integer then $|\hat f|^{p'} = \hat F$ where $\|F\|_1 \le
\|f\|_1^{p'}$. Using the same proof as for $p' = 2$ we get for $d \ge 4$
$$\|\hat f\|_{p'}^{p'} \le C(\int\limits_{\rho \in SO(d)}\|\hat
g_{\rho}\|_{p'}^{p'} d\rho +
\|f\|_1^{p'})$$
and
for $d \ge 5$
$$\int\limits_{\rho \in SO(d)}\|\hat g_{\rho}\|_{p'}^{p'} d\rho  \le C(\|\hat
f\|_{p'}^{p'} +
\|f\|_1^{p'}).$$

~

{\bf 2. Case $1 \le p < \frac {2d}{d + 2}$.}\\

We will generalize {\bf Theorems 1 and 1}$'$ in the spirit of the Stein-Tomas
Theorem
(\cite{DC}, Chapter 6.5).

~

{\bf Theorem 2}:
{\it let $d \ge 4$ and let $f \in L^1({\Bbb {R}} ^d) \cap L^p({\Bbb {R}} ^d)$
where $1 \le p < \frac
{2d} {d + 2}$. If periodizations of f
\begin{eqnarray*} g_{\rho}(x) = \sum \limits _{\nu \in {\Bbb {Z}} ^d} f(\rho (x
- \nu))\end{eqnarray*}
are in $L^2([0,1]^d)$ for almost all rotations $\rho \in SO(d)$ and
\begin{eqnarray*} G^2 = \int\limits_{\rho \in SO(d)}\|g_{\rho}\|_2^2d\rho <
\infty
\end{eqnarray*}
then $f \in L^2({\Bbb {R}} ^d)$:
\begin{eqnarray}\|f\|_2 \le C(G + \|f\|_p)\label{G78}\end{eqnarray}
where $C$ depends only on $d$ and $p$.}

~

It will follow from the proof (see \bref{g21}) that we can replace
$\int\limits_{\rho \in SO(d)}\|g_{\rho}\|_2^2 d\rho$ with $\int\limits_{\rho \in
SO(d)}\|g_{\rho} - \hat g (0)\|_2^2 d\rho$ in {\bf Theorem 2}, which is the norm
of $g$ in the quotient space
$L^2([0,1]^d \times SO(d))$ modulo constants.
We will also obtain an inverse theorem.

~

{\bf Theorem 2}$'$:
{\it let $d \ge 5$ and let $f \in L^1({\Bbb {R}} ^d) \cap L^2({\Bbb {R}} ^d)$
and
$1 \le p < \frac {2d} {d + 2}$ and let $g_{\rho}$ be periodizations of f
\begin{eqnarray*} g_{\rho}(x) = \sum \limits _{\nu \in {\Bbb {Z}} ^d} f(\rho (x
- \nu))\end{eqnarray*}
then $g_{\rho} \in L^2([0,1]^d)$ for almost all rotations $\rho \in SO(d)$ and
\begin{eqnarray}\int\limits_{\rho \in SO(d)}\|g_{\rho} - \hat g_{\rho} (0)\|_2^2
d\rho  \le C(\|f\|_2 + \|f\|_p)^2\label{G79}\end{eqnarray}
where $C$ depends only on $d$ and $p$.}\\

Since Schwartz functions are dense in $L^p({\Bbb {R}} ^d) \cap L^2({\Bbb {R}}
^d)$ it
follows from {\bf Theorem 2}$'$ that we can define periodizations $g_{\rho}$ of
$f \in L^p({\Bbb {R}} ^d) \cap L^2({\Bbb {R}} ^d)$ where
$1 \le p < \frac {2d} {d + 2}$ for a.e. $\rho \in SO(d)$ as elements of the
quotient space of $L^2([0,1]^d)$ modulo constants. \\

We say that $f \in L^p({\Bbb {R}} ^d)$ has periodizations $g$ in the quotient
space
$L^2([0,1]^d \times SO(d))$ modulo constants if there exists a sequence of
Shwartz functions $f_k$ converging to $f$ in $L^p({\Bbb {R}} ^d)$ and such that
$g_k \rightarrow g$ in the quotient space
$L^2([0,1]^d \times SO(d))$ modulo constants. From {\bf Theorem 2} we conclude
that $f \in L^2({\Bbb {R}} ^d)$ and $f_k \rightarrow f$ in $L^2({\Bbb {R}} ^d)$.
It follows from {\bf Theorem 2}$'$ that $g$ is a well-defined element of the
quotient space
$L^2([0,1]^d \times SO(d))$ modulo constants.\\

{\bf Remarks}: 1. As the following example shows, we can not replace
$\int\limits_{\rho \in SO(d)}\|g_{\rho} - \hat g (0)\|_2^2 d\rho$ with
$\int\limits_{\rho \in SO(d)}\|g_{\rho}\|_2^2 d\rho$ in {\bf Theorem 2}$'$ when
$p > 1$. Let $\phi: \r \rightarrow {\Bbb {C}}$ be a Schwartz function supported
in $B(0,1)$ such that $\phi(0) = 1$. Put $\hat f(x) =
\phi(\frac x {\epsilon})$. Then $$g_{\rho} = \hat f (0) = 1$$
 but
$$\|f\|_p = \epsilon ^ {\frac d {p'}} \|\check \phi\|_p.$$

2. The next example from (\cite{DC}, Chapter 6.3) shows that $p$ can not be
greater
than $\frac {2d + 2} {d + 3}$ in {\bf Theorem 2}$'$. Put
\begin{eqnarray}\hat f(x_1,...,x_d) = \phi(\frac {x_1 - 1} {\epsilon ^2}, \frac
{x_2} {\epsilon},..., \frac {x_d} {\epsilon}) \label{EX}\end{eqnarray}
where $\phi: \r \rightarrow {\Bbb {C}}$ is a Schwartz function supported in
$B(0,2)$ such that $\phi = 1$ in $B(0,1)$. Then
\begin{eqnarray*} \int\limits_{\rho \in SO(d)}\|g_{\rho} - \hat g (0)\|_2^2
d\rho &=& 2d \int_{|\xi| = 1} |\hat f(\xi)|^2 d\sigma (\xi)\\
 &\ge& C \epsilon ^ {d-1}\end{eqnarray*} but
$$\|f\|_p^2 = \epsilon ^ {\frac {2d + 2} {p'}}\|\check \phi\|_p.$$
It is an open question whether {\bf Theorems 2} and {\bf 2}$'$ are valid when
$\frac
{2d}{d + 2} \le p < \frac {2d + 2}{d + 3}$. We discuss this further in {\bf
Remark 2} at the
end of the paper.

~

Proof of {\bf Theorem 2}: \\
The proof is quite similarly to that of {\bf Theorem 1}. We will replace {\bf Lemma 1}
with

~

{\bf Lemma 3} {\it
Let $q: {\Bbb {R}} \rightarrow {\Bbb {R}}$ be a Schwartz function supported in
$[\frac 12, 2]$, let $f \in L^p(\r)$ where $1 \le p \le 2$ and let $b \in
[0,1)$.
Define $H_1:{\Bbb {R}} \rightarrow {\Bbb {C}}$
\begin{eqnarray*} H_1 (t) = \frac 1 {\w} h_1(\w) q(\frac {\w} N).\end{eqnarray*}
Then for large enough $N$ we have

\begin{eqnarray} \sum \limits _{\nu \neq 0} |\hat H_1 (\nu)| \le \frac
{C\|f\|_p^2} {N} \label{H111}\end{eqnarray} where $C$ depends only on $q$ and $d$.}

~

Proof of {\bf Lemma 3}:\\

The only difference in the proof is how to obtain an inequality analogous to \bref{dif}.
Using Young's inequality we have $\|f* \tilde f\|_q \le \|f\|_p^2$ where $1 + \frac 1 q =
\frac 2 p$. Therefore $|\int (f * \tilde f)(x) w(x) dx| \le \|f\|_p^2
 \|w\|_{q'}$. Substituting derivatives of $\widehat {d \sigma _t} (x)
\chi_{\{|x| \le 1\}}$ with respect to $t$
instead of $w$, we get the desired inequality
\begin{eqnarray}|h_1^{(k)} (t)| \le Ct^{d-1}\|f\|_p^2\end{eqnarray}
where $t \ge 1$ and $C$ depends only on $k$ and $d$. \hfill$\square$\\

The main difficulty is to prove a lemma analogous to {\bf Lemma 2}:

~

{\bf Lemma 4} {\it
Let $q: {\Bbb {R}} \rightarrow {\Bbb {R}}$ be a Schwartz function supported in
$[\frac 12, 2]$, let $f \in L^p({\Bbb {R}} ^d)$ where $1 \le p < \frac {2d} {d +
2}$
 and let $b \in [0,1)$. Define $H_{2,N}:{\Bbb {R}} \rightarrow {\Bbb {C}}$
\begin{eqnarray*} H_{2, N} (t) = \frac 1 {\w} h_2(\w) q(\frac {\w}
N).\end{eqnarray*}
Then we have
\begin{eqnarray} \sum \limits _{\nu \neq 0}|\sum \limits _{j \ge 0} \hat H_{2,
2^j} (\nu)|
 \le C\|f\|_p^2\label{H222}\end{eqnarray}
with $C$ depending only on $p$, $q$ and $d$.}

~

Proof of {\bf Lemma 4}:\\

Recall from \bref{H3} that

\begin{eqnarray*} \hat H_{2, N} (\nu) = 2\int (f* \tilde f)(x)D_{N, \nu}(x) dx
\end{eqnarray*}
where \begin{eqnarray} D_{N, \nu}(x) = \chi_{\{|x| > 1\}}e^{i2\pi \nu b}\int
Nq(t)
e^{-i2\pi \nu (Nt)^2} (Nt)^{d-1}\widehat
{d\sigma}(Ntx)dt.\label{D99}\end{eqnarray}
Denote by
\begin{eqnarray} K_{\nu}(x) = \sum \limits _{l \ge 0}D_{2^l, \nu}(x)
\label{K99}.
\end{eqnarray}
Then \begin{eqnarray} |\sum \limits _{l \ge 0} \hat H_{2, 2^l} (\nu)| &=& 2|\int
(f* \tilde f)(x) \sum \limits _{l \ge 0}D_{2^l, \nu}(x) dx| \nonumber \\
&=& 2|\int \tilde f(x) (K_{\nu}*f)(x) dx| \nonumber \\
&\le& 2\|f\|_p \|K_{\nu}*f\|_{p'}\label{H20}.\end{eqnarray}
If $p' = \infty$ or $p' = 2$ we have
\begin{eqnarray*} \|K_{\nu}*f\|_{\infty} \le \|K_{\nu}\|_{\infty}\|f\|_1 \\
\|K_{\nu}*f\|_{2} \le \|\hat K_{\nu}\|_{\infty}\|f\|_2.\end{eqnarray*}
First we will show that
\begin{eqnarray} \|K_{\nu}\|_{\infty} &\le& \|\sum \limits _{l \ge 0}
|D_{2^l, \nu}|(x)\|_{\infty} \nonumber \\
&\le& C\nu ^{-\frac{d}2}\label{K1}.\end{eqnarray}
It follows from \bref{phi1} that
\begin{eqnarray} |D_{N, \nu}(x)| \le \frac{N^{\frac{d+1} 2}}{|x|^{\frac{d-1}
2}}\begin{cases} C(\nu N^2)^{-\frac 12}&\text{if $N \in [\frac {|x|}
{4\nu},\frac {|x|} {\nu}]$}\\
C_k(\nu N^2)^{-k-\frac 12}&\text{if $N \notin [\frac {|x|} {4\nu},\frac {|x|}
{\nu}]$}
\end{cases}
.\label{D1}\end{eqnarray}
If $\nu > 0$ then the number of diadic $N \in [\frac {|x|} {4\nu},\frac {|x|}
{\nu}]$ is
 at most $3$. If $\nu < 0$ then there are no $N$ in $[\frac {|x|} {4\nu},\frac
{|x|} {\nu}]$. Therefore choosing $k \ge \frac {d-1}2$ and summing \bref{D1}
over all diadic $N$ we have
\begin{eqnarray*}\sum \limits _{l \ge 0}|D_{2^l, \nu}(x)| \le C|\nu|^{-\frac {d}
2}
\end{eqnarray*}
with $C$ depending only on $d$ and $q$.\\

Now we will show that
\begin{eqnarray} \|\hat K_{\nu}\|_{\infty} \le \|\sum \limits _{l \ge 0}|\hat
D_{2^l, \nu}|(y)\|_{\infty} \le C \label{K2}.\end{eqnarray}
Since supp $\phi \in [\frac 12, 2]$ we can re-write \bref{phi} for a stronger
version of
the method of stationary phase (\cite{H}, Theorems 7.6.4, 7.6.5, 7.7.3)
\begin{eqnarray}|\int \phi(t)e^{-i2\pi \nu N^2 (t-\frac {|x|}{2\nu N})^2}dt -
\sum \limits _{j=0}^{k-1}c_j(\nu N^2)^{-j-\frac 12}\phi^{(2j)}(\frac {|x|}{2\nu
N})| \le \frac
{c_k(\nu N^2)^{-k-\frac 12}} {max(1, \frac {|x|} {8N\nu})^k}
\nonumber\end{eqnarray}
where $c_j$ are some constants.
Therefore
\begin{eqnarray}D_{N, \nu}(x) = \chi_{\{|x| > 1\}}\frac {N^{\frac {d+1} 2}}
{|x|^{\frac {d-1} 2}}e^{i2\pi \frac {|x|^2}{4\nu}}\sum \limits
_{j=0}^{k-1}c_j(\nu N^2)^{-j-\frac 12}\phi^{(2j)}(\frac {|x|}{2\nu N}) +
\phi_k(x)\label{D10}\end{eqnarray}
where $|\phi_k(x)| \le \chi_{\{|x| > 1\}}\frac {N^{\frac {d+1} 2}} {|x|^{\frac
{d-1} 2}}\frac {c_k(\nu N^2)^{-k-\frac 12}} {max(1, \frac {|x|} {8N\nu})^k}$.
Choosing $k \ge \frac {d + 2} 2$ we have
\begin{eqnarray}\|\hat \phi_k\|_{\infty} &\le& \|\phi_k\|_1\nonumber \\
&=& \int\limits_{|x| \le 8\nu N}|\phi_k| dx + \int\limits_{|x| > 8\nu N}|\phi_k|
dx \nonumber \\
&\le& \frac C N \label{phik}\end{eqnarray}
where $C$ depends only on $d$ and $q$. We can assume that $\nu > 0$ since $D_{N,
\nu}(x) = \phi_k(x)$ for $\nu < 0$. We can also ignore $\chi_{\{|x| > 1\}}$ in
front of the sum in
\bref{D10} because if $\frac {|x|}{2\nu N} \in [\frac 12, 2]$, then $|x| \ge \nu
N \ge 1$.
 We will consider only the zero term in the sum. The other terms can be treated
similarly.
The Fourier transform of
\begin{eqnarray*}\frac {N^{\frac {d+1} 2}} {|x|^{\frac {d-1} 2}}e^{i2\pi \frac
{|x|^2}{4\nu}}
(\nu N^2)^{-\frac 12}\phi(\frac {|x|}{2\nu N})\end{eqnarray*}
at point $y$ is equal to
\begin{eqnarray}N^{\frac {d+1} 2}(2\nu N)^{\frac {d+1} 2}(\nu N^2)^{-\frac
12}\int_{\r} \psi(|x|)e^{i2\pi \nu N^2 |x|^2} e^{-i2\pi 2\nu N x \cdot y} dx =
\nonumber \\
C(\nu N^2)^{\frac d 2}e^{-i2\pi \nu|y|^2}\int_{\r} \psi(|x|)e^{i2\pi \nu N^2 |x - \frac y
N|^2}dx\label{psi10}\end{eqnarray} where $\psi(t) = \phi(t, 2\nu N t) t^{-\frac{d-1}2}$
is a Schwartz function supported in $[\frac 12, 2]$ whose derivatives and the function
itself are bounded uniformly in $t$, $\nu$ and $N$ (see the remarks after \bref{a1}). The
same is true about partial derivatives of $\psi(|x|)$. Applying the stationary phase
method for $\r$ (\cite{H}, Theorem 7.7.3) we get
\begin{eqnarray}|\int_{\r} \psi(|x|)e^{i2\pi \nu N^2 |x - \frac y N|^2}dx| \le
\begin{cases} C(\nu N^2)^{-\frac d 2}&\text{if $N \in [\frac {|y|} {2},2|y|]$}\\
C_k(\nu N^2)^{-k-\frac d 2}&\text{if $N \notin [\frac {|y|}
{2},2|y|]$}\end{cases}
.\label{psi1}\end{eqnarray}
Therefore the absolute value of \bref{psi10} can be bounded from above by:
\begin{eqnarray} \le \begin{cases} C&\text{if $N \in [\frac {|y|} {2},2|y|]$}\\
C_k(\nu N^2)^{-k}&\text{if $N \notin [\frac {|y|} {2},2|y|]$}\end{cases}
.\label{psi2}\end{eqnarray}
Similar inequalities hold for Fourier transforms for the rest of the terms in
the sum in \bref{D10}.
The number of diadic $N \in [\frac {|y|} {2},2|y|]$ is bounded by $3$. Using
\bref{phik},
 choosing $k \ge 1$ in \bref{psi2} and summing over all diadic $N$ we get
\begin{eqnarray}\sum \limits_{l \ge 0}|\hat D_{2^l, \nu}(y)| \le C
\label{D5}\end{eqnarray}
with $C$ depending only on $d$ and $q$.
Using \bref{K1} and \bref{K2} and interpolating between $p=1$ and $p=2$, we
obtain
\begin{eqnarray*} \|K_{\nu}*f\|_{p'} \le C\nu ^{-\alpha
_p}\|f\|_p\end{eqnarray*} where $\alpha_p = \frac d 2 \frac {2-p}p$. We have $\alpha_p >
1$ if $p < \frac {2d}{d + 2}$. Summing \bref{H20} over all $\nu \neq 0$, we get the
desired inequality
\begin{eqnarray*} \sum \limits _{\nu \neq 0}|\sum \limits _{j \ge 0} \hat H_{2,
2^j} (\nu)|
 \le C\|f\|_p^2.\end{eqnarray*}\hfill$\square$

~

Now we are in a position to proceed with the proof of {\bf Theorem 2}. The proof
is almost
 the same as the one of
{\bf Theorem 1}.  We also need to replace inequality \bref{h16} with the inequality
\begin{eqnarray*}\int\limits_0^1 h(t)dt &=& \int\limits_{|y| \le 1} |\hat f
(y)|^2 dy \\
&\le& C\|\hat f\|_{p'}^2 \\
&\le& C\|f\|_p^2\end{eqnarray*}
where $p \le 2$ and $C$ depends only on $d$. An argument similar to the one used
to get
 \bref{h35}, \bref{h15} and \bref{f10} yields the desired inequality
\begin{eqnarray*}\int |\hat f(x)|^2 dx &=& \int\limits_0^{\infty} |\hat
f(\xi)|^2 d\sigma_{t}(\xi)dt\\
&=& \int\limits_0^{\infty} h(t) dt \\
&\le& C(G^2 + \|f\|_p^2)
\end{eqnarray*}
with $C$ depending on $d$ and $p$. Note that the interchange of summation by $\nu$
 and $N$ is not a problem. \hfill$\square$\\

~

The proof of {\bf Theorem 2}$'$ is the same (see the argument before {\bf Theorem 1}$'$).
The important thing is that now we exclude $\hat g_{\rho} (0) = \hat
f(0)$ in \bref{g20} now.\\

{\bf Final remarks}: 1. We can further generalize {\bf Theorem 2}$'$. Fix some $q \in [1,
\frac {2d}{d+2})$. Applying complex interpolation between the trivial $p=1$ and $p=2$, we
obtain the
 following result for $1 < p < 2$: let $d
\ge 5$ and
 let $f \in L^{\frac {qp} 2}({\Bbb {R}} ^d) \cap L^p({\Bbb {R}} ^d)$ and let
$g_{\rho}$ be
periodizations
 of f
\begin{eqnarray} g_{\rho}(x) = \sum \limits _{\nu \in {\Bbb {Z}} ^d} f(\rho (x
- \nu))\label{g50}\end{eqnarray}
then $g_{\rho} \in L^p([0,1]^d)$ for almost all rotations $\rho \in SO(d)$ and
$$\left (\int\limits_{\rho \in SO(d)}\|g_{\rho}-\hat g_{\rho}(0)\|_p^{p'} d\rho
\right) ^ {\frac 1 {p'}} \le
C\|f\|_p^{2 - p}(\|f\|_p + \|f\|_{\frac {qp} 2})^{p - 1}$$
where $C$ depends only on $d$ and $q$. Choosing $q = \frac 2 p$ in the above
inequality we obtain the following generalization of {\bf Theorem 1}$'$: if
$\frac {d + 2}{d} < p \le 2$ and $d \ge 5$ then
$$\left (\int\limits_{\rho \in SO(d)}\|g_{\rho}\|_p^{p'} d\rho \right) ^ {\frac
1 {p'}} \le
C_p(\|f\|_p + \|f\|_1).$$\\

If $p'$ is an even integer then $|\hat f|^{p'} = \widehat {F * \tilde F}$ where $\|F\|_r
\le \|f\|_q^{\frac {p'} 2}$ with $\frac {p'} 2 - 1 + \frac 1 r = \frac {p'/2} q$.
Repeating the same arguments as for $p' = 2$ we obtain for $d \ge 4$

$$\|\hat f\|_{p'}^{p'} \le C (\|f\|_q^{p'} + \int\limits_{\rho \in SO(d)}\|\hat
g_{\rho}\|_{p'}^{p'} d\rho)$$
and for $d \ge 5$
$$\int\limits_{\rho \in SO(d)}\|\hat g_{\rho} - \hat g_{\rho}(0)\|_{p'}^{p'}
d\rho  \le C (\|f\|_q^{p'} + \|\hat f\|_{p'}^{p'})$$
where $1 \le q < \frac {\frac {p'} 2}{\frac {p'} 2 - 1 + \frac {d + 2} {2d}}$
and $C$ depends only on $d$ and $q$.

2. Conditionally on the exponent pair conjecture (\cite{M}, Chapter 4,
Conjecture 2) we can
 clarify what happens when $ \frac {2d}{d+2} \le p < \frac{2d + 2}{d+3}$. In our
case the conjecture says that
\begin{eqnarray}|\sum\limits_{n \le \nu \le m} e^{i\frac {|x|^2} {\nu}}| \le
C_{\epsilon} |x|^{\epsilon} n^{\frac 12} \label{EP}\end{eqnarray}
where $m \le 2n$ and $|x| \ge n$. Let $\beta(x) = \max(1, |x|)$.\\

{\bf Proposition 1.} {\it {\bf Theorems 2 and 2}$'$ hold if we replace $\|f\|_p$ with
$\|\beta ^ {\epsilon}f\|_p$ and if $p < \frac{2d + 2}{d+3}$, provided the
conjecture \bref{EP} is valid.}\\

Using the example \bref{EX} we can show that the {\bf Proposition 1} is sharp up to
$\epsilon$
 in the range of $p$ for the estimate \bref{G79}.\\

{\bf Proof of Proposition 1}:\\
The main issue is to improve the result of {\bf Lemma 4}. Denote by
\begin{eqnarray}L_{j}(x) = \sum\limits_{2^j \le \nu \le
2^{j+1}}K_{\nu}(x)\label{L}.
\end{eqnarray}
Using summation by parts we obtain from \bref{D10}, \bref{K99} and \bref{EP}
that
\begin{eqnarray*}|L_{j}(x)| &=& |\sum\limits_{2^j \le \nu \le
2^{j+1}}K_{\nu}(x)| \\
 &\le& C_{\epsilon}|x|^{ \epsilon}2^{-\frac {d-1} 2 j}.\end{eqnarray*}
We will deal with the following expression instead of \bref{H20}
\begin{eqnarray*}|\int \tilde f(x) (L_{j}*f)(x) dx| &=&
 |\int \tilde f(x)\beta^{\epsilon} \frac {(L_{j}*f)(x)} {\beta^{\epsilon}} dx|
\\
&\le& \|\beta^{\epsilon}f(x)\|_p \|\frac {(L_{j}*f)}
{\beta^{\epsilon}}\|_{p'}.\end{eqnarray*}
If $p' = \infty$ or $p' = 2$ we have
\begin{eqnarray*} \|\frac {(L_{j}*f)} {\beta^{\epsilon}}\|_{\infty} &\le&
\|\frac {\int_{|y|
\le |x|}|L_{j}(x-y)|\cdot |f(y)|dy + \int_{|y| \ge |x|}|L_{j}(x-y)|\cdot
|f(y)|dy} {\beta^{\epsilon}(x)}\|_{\infty}\\
&\le& \|\frac {L_j} {\beta^{\epsilon}}\|_{\infty}\|\beta^{\epsilon}f\|_1 \\
&\le& C_{\epsilon}2^{-\frac {d-1} 2 j}\|\beta^{\epsilon}f\|_1, \\
\|\frac {(L_{j}*f)} {\beta^{\epsilon}}\|_{2} &\le& \|L_{j}*f\|_{2}\\
 &\le& \|\hat L_{j}\|_{\infty}\|f\|_2\\
&\le& C 2^j \|f\|_2.\end{eqnarray*}
Interpolating between $p=1$ and $p=2$, we obtain
\begin{eqnarray*} \|\frac {(L_{j}*f)} {\beta^{\epsilon}}\|_{p'} \le
C_{\epsilon}2^{-j\alpha _p}\|\beta^{\epsilon}f\|_p\end{eqnarray*} where $\alpha_p = \frac
{d + 1} p -  \frac {d+3} 2$. We have $\alpha_p > 0$ if $p <
\frac {2d + 2}{d + 3}$.\hfill$\square$ \\

3. Concerning the lower dimensional cases we can use the following results from
the Number
Theory:
$$r_3(n) \le C n^{\frac 12}\ln n \ln \ln n,$$
$$r_4(n) \le C n \ln \ln n.$$
See for example (\cite{B}).
There is an infinite arithmetic progression, e.g. $n = 8k + 1$, such that
$$r_3(n) \ge C_{\epsilon} n^{\frac 12 - \epsilon}.$$
See for example (\cite{GCC}).
Then {\bf Theorem 2} holds when $d= 3$ and {\bf Theorem 2$'$} holds when $d= 3$
or $d=4$ if we replace $$\|g_{\rho}\|_2^2 = \sum \limits _{m \in \z} |\hat
g_{\rho}(m)|_2^2$$ with
$$\sum \limits _{m \in \z} |m|^{\epsilon}|\hat g_{\rho}(m)|_2^2,$$
$$\sum \limits _{m \in \z, |m| > 3}  \frac {|\hat g_{\rho}(m)|_2^2}{\ln |m| \ln
\ln |m|}$$ or
$$\sum \limits _{m \in \z, |m| > 3}\frac {|\hat g_{\rho}(m)|_2^2}{\ln \ln |m|}$$
correspondingly.
\\
Using a technique similar to the one in the proof of {\bf Proposition 1}, second remark
we also obtain the following results in lower dimensions:

$$\|\beta \cdot \hat f\|_2 \lesssim \|f\|_p + \|g\|_2 \mbox{ when $d = 3$ and $\beta (x) = \min (1, \frac 1 {|x|^{\epsilon}})$},$$

$$\|g - \hat g (0)\|_2 \lesssim \|\beta \cdot f\|_p + \|\beta \cdot \hat f\|_2
\mbox{ when $d = 3$ and $\beta (x) = \max (1, \ln_+ |x| \ln_+ (\ln_+ |x|))$},$$

$$\|g - \hat g (0)\|_2 \lesssim \|\beta \cdot f\|_p + \|\beta \cdot \hat f\|_2
\mbox{ when $d = 4$ and $\beta (x) = \max (1, \ln_+ (\ln_+ |x|))$},$$
and $1 \le p < \frac {2d}{d + 2}$ in all three cases.
~

\begin{center}
{\bf Acknowledgements}
\end{center}
This work was inspired by useful discussions with Thomas Wolff.

~

\end{document}